\documentstyle[11pt]{article}

\newtheorem{theorem}{Theorem}
\newtheorem{definition}[theorem]{Definition}
\newtheorem{lemma}[theorem]{Lemma}
\newtheorem{corollary}[theorem]{Corollary}
\newtheorem{proposition}[theorem]{Proposition}
\newtheorem{example}[theorem]{Example}

\newcommand\ad{\mathop{\it ad}\nolimits}

\newcommand\Lie{\mathop{\sl Lie}\nolimits}

\newenvironment{proof}%
{\rm \trivlist \item[\hskip\labelsep{\sc Proof.}]}%
{\hspace*{\fill}$\Box$ \endtrivlist}

\newcommand \frak{\bf \sf} 
\newcommand \Bbb {\bf }    

\newcommand \fg {{\frak g}}
\newcommand \fn {{\frak n}}
\newcommand \fm {{\frak m}}
\newcommand \fh {{\frak h}}
\newcommand \fs {{\frak s}}
\newcommand \fl {{\frak l}}

{\rm \trivlist \item[\hskip\labelsep{\sc Sketch of Proof.}]}%
{\hspace*{\fill}$\Box$\endtrivlist}

\newcommand{\lefthook}{\mbox{$\rule{8pt}{.5pt}\rule{.5pt}{6pt}\, $}}
\begin{document}

\title{\rm Proper group actions and  symplectic stratified spaces }

\author{L. Bates\thanks{Department of Mathematics, University of
Calgary, Calgary, Alberta, Canada T2N 1N4.} and E.
Lerman\thanks{Department of Mathematics, Massachusetts Institute of
Technology, Cambridge, Ma 02139, USA.}
\thanks {Supported by an NSF postdoctoral fellowship}}
\date{}
\maketitle
\begin{abstract}
Let $(M,\omega)$ be a Hamiltonian $G$-space with a momentum map $F:M
\to {\frak g}^*$. It is well-known that if $\alpha$ is a regular
value of $F$ and $G$ acts freely and properly on the level set
$F^{-1}(G\cdot \alpha)$, then the reduced space $M_\alpha
:=F^{-1}(G\cdot \alpha)/G$ is a symplectic manifold. We show
 that if
the regularity assumptions are dropped the space $M_\alpha $ is a
union of symplectic manifolds, and that the symplectic manifolds fit
together in a nice way.  In other words the reduced space is a {\em
symplectic stratified space}. This extends results known for the
Hamiltonian action of compact groups.
\end{abstract}

\section*{\rm{Introduction}}
Reduction of the number of degrees of freedom of a symmetric Hamiltonian
system has a long history.  The modern formulation of reduction is due
to Meyer \cite{My} and to Marsden and Weinstein~\cite{MW}.
We recall their result.  One starts with a symplectic manifold
$(M, \omega)$, a Hamiltonian action of a Lie group $G$ and a corresponding
equivariant momentum map $F: M \to {\fg}^*$.  Let ${\cal O}$ be a coadjoint
orbit of $G$.  If the momentum map is transversal to the orbit,
then the preimage $F^{-1}({\cal O})$ of the orbit is a submanifold  of
$M$ and the action of the Lie group $G$ on the preimage is locally free.
Assume that this action is actually free and that the orbit map
$F^{-1}({\cal O} )\to F^{-1}({\cal O})/G$ is a fibration.
The reduction theorem
says that the orbit space $M_{\cal O} := F^{-1}({\cal O})/G$ is a symplectic
manifold.  The restriction of  a smooth $G$ invariant function $h$ on $M$ to
 the preimage of the orbit descends to a smooth function $h_{\cal O}$ on
the reduced space $M_{\cal O}$.  Moreover, the Hamiltonian flow of $h$ on
$F^{-1}({\cal O})$ descends to the Hamiltonian flow of $h_{\cal O}$ on the
reduced space.

It turns out that often the action is only locally free, so at best
the reduced spaces are symplectic orbifolds.  This already suggests that
the category of symplectic manifolds is too restrictive for
Hamiltonian dynamics.  More generally one would like to get rid of the
the transversality hypothesis in the reduction procedure.  One reason for
this desire is that the more symmetry the point of a system has
the more singular the momentum map is at this point.  Of course,
the symmetric points are not generic, but they are very important
in understanding of the dynamics of the system.  Another reason
is that one would like to understand the change in the topology of the
reduced space as one crosses the critical values of the momentum map.

For a number of years the reduction at singular values of the momentum
map has been problematic. In 1981 Arms, Marsden and Moncrief \cite{AMM}
showed that under some assumptions the {\em set} $F^{-1}({\cal O})/G$ is
a union of symplectic manifolds and that the flow of invariant
Hamiltonians on the level set $F^{-1}({\cal O})$ of the momentum map
descends to the flow of the reduced Hamiltonians on these symplectic
manifolds.  Yet this observation didn't gain use and is not well known.
In fact it has been rediscovered at least once \cite{Otto}.
Many reduction schemes have been proposed since 1981. A number of
them are compared in \cite{AGJ}.

Our approach to  reduction is the one proposed in \cite{SL} and
\cite{lms}.  Namely, for a point $\alpha $ in the dual of Lie algebra
of $G$, the reduced space at $\alpha$ {\em is} the topological space
$M_\alpha = M_{G\cdot \alpha} := F^{-1}(G\cdot \alpha )/G$, where
$G\cdot \alpha $ is the coadjoint orbit through $\alpha$.  In general
this topological space can be quite horrible, as we shall see shortly.
One of the main points of the paper is that we only need to make two
assumptions --- that the action is proper and that the coadjoint
orbits of our group are locally closed --- to guarantee that the
reduced spaces are manageable.  By `manageable' we mean that
Hamilton's equations hold and the geometry of the reduced space is
reflected in the dynamics.

We will also show that in analogy with symplectic orbifolds (which are
modeled on a symplectic vector space modulo a finite group) our
reduced spaces are modeled on symplectic vector spaces reduced at zero
with respect to a linear action of a compact group.  This extends the
results of \cite{SL} and \cite{CushSj} which proved the above
assertions for the case of the compact symmetry group.  One motivation for
the extension is to push the methods of \cite{SL} as far as they would go.
Another motivation
for this extension comes from field theory, where the symmetry groups
are not compact. Yet some field theories such as Yang-Mills in bounded
domains do not satisfy the assumptions of the Arms-Marsden-Moncrief
theory (for which field theory appears to be a primary motivation),
but the gauge group still acts properly, and a large portion of the
finite dimensional results can still be established \cite{BSS}.

We now  briefly describe the organization of the paper.
\begin{enumerate}
\item We start out by defining an algebra of ``smooth functions'' on the
reduced space with a natural Poisson bracket.
The bracket allows us to define Hamiltonian flows of smooth functions on
the reduced space.  If the smooth functions on the reduced space separate
points the flows are unique.

\item The Hamiltonian flows of smooth functions  preserve
the decomposition of  the reduced space induced by the orbit type decomposition
of the original manifold.

\item  Local normal form computations show that
\begin{enumerate}
\item the orbit type decomposition of the reduced space is a decomposition into
symplectic manifolds;

\item  the embeddings of these manifolds ({\em the symplectic pieces}) into
the reduced space are Poisson maps;

\item  the group generated by the Hamiltonian flows of functions on the reduced
space acts transitively on the connected components of the symplectic pieces;

\item  consequently, the Poisson algebra of smooth functions on the reduced
space carries all the information about the decomposition of the reduced space
into symplectic pieces.
\end{enumerate}

\item  The last fact allows us to define isomorphisms of reduced spaces
in terms of the corresponding isomorphisms of Poisson algebras of functions.
We can also define local isomorphisms.

\item  Local normal form computations show that a reduced space is locally
isomorphic to a symplectic vector space reduced at zero with respect to a
linear action of a compact group.  This symplectic vector space is the
maximal symplectic subspace of the slice to the corresponding orbit in
the original manifold.

\item  It follows that the decomposition of the  reduced space by orbit type
is a stratification and that the local structure of a stratification can be
read off from the slice representation.

\item We use the local normal form computation to show that the strata
of the reduced space can individually be obtained by Marsden-Weinstein-Meyer
reduction.  This provides us with a way to reconstruct the original dynamics
from the dynamics on the reduced space.

\item We conclude by showing how one can use symplectic cross-sections to
factor out the coadjoint orbit directions.

\end{enumerate}

\section{Dynamics on the reduced space}

Consider a symplectic manifold $M$ with a Hamiltonian action of a Lie
group $G$ and let $F:M \to {\frak g}^*$ be a corresponding equivariant
momentum map.  Fix a coadjoint orbit ${\cal O}$ of $G$.  We {\em
define} the corresponding reduced space $M_{\cal O}$ to be the
topological quotient of the subset $F^{-1}({\cal O})$ of $M$ by the
action of the group $G$,
$$
M_{\cal O}:= F^{-1}({\cal O})/G.
$$
We have not made enough assumptions to guarantee that the set
$F^{-1}({\cal O})$ is a manifold or that the quotient space
$M_{\cal O}$ is nice.

\begin{example} {\rm
Consider an irrational flow on a torus ${\Bbb R}\times {\Bbb T}^2 \to {\Bbb
T}^2 $
generated by a vector $\xi $ in the Lie algebra of ${\Bbb T}^2$.  The flow
lifts to a
Hamiltonian action on the cotangent bundle of the two torus.  The reduced space
at zero
$M_0$ is homeomorphic to ${\Bbb T}^2 /{\Bbb R}\times {\Bbb R} \xi ^\circ $
where
${\Bbb R} \xi ^\circ $ is the annihilator in the dual of the Lie algebra of the
torus of
the line through $\xi$.  The reduced space is not Hausdorff.

We note for future reference that the space of  functions on the cotangent
bundle of the
two torus that are invariant under the flow is isomorphic to the space of
functions
on ${\Bbb R}^2$ and that the Poisson bracket of two invariant functions is
zero.
}
\end{example}

Our first step in defining the dynamics on the reduced space (in this
we are following
\cite{ACG}) is to define a Poisson algebra of ``smooth functions'' on the
reduced space.
Since the restriction of a smooth invariant function on the manifold $M$ to
 the set $F^{-1}({\cal O})$ descends to a continuous function on the quotient
$F^{-1}({\cal O})/G = M_{\cal O}$, we {\em define} the smooth functions on the
reduced space to be
these restrictions,
\begin{eqnarray*}
	C^\infty (M_{\cal O}) &:= &\left. C^\infty (M)^G \right| _{F^{-1}({\cal O})}\\
			      &:= & C^\infty (M)^G /{\cal I},
\end{eqnarray*}
Here $C^\infty(M)^G$ is the algebra of smooth $G$-invariant functions
on the manifold $M$, and ${\cal I = I}(F^{-1}({\cal O}))$ is the ideal
of invariant functions that vanish on the set $F^{-1}({\cal O})$.

To show that the algebra of smooth functions $C^\infty (M_{\cal O})$
is a Poisson algebra, we need to check that ${\cal I}$ is not only an
ideal under multiplication of functions but also an ideal with respect
to the Poisson bracket (recall that the $G$ invariant functions form a
Poisson subalgebra of $C^\infty (M)$).  The fact that ${\cal I}$ is a Poisson
ideal follows from  Lemma~\ref{lemma1} below.

\begin{lemma}\label{lemma1}
Let $M$ be a symplectic (or, more generally, a Poisson) manifold,
and ${\cal A}$ be a Poisson subalgebra of $C^\infty (M)$.  Suppose that the
Hamiltonian flows of functions in ${\cal A}$ preserve a subset $X$ of the
manifold $M$.  Then the ideal
$$
	{\cal I}(X) := \{f \in {\cal A} : f|_X = 0\}
$$
of functions in ${\cal A}$ that vanish on $X$ is a Poisson ideal of ${\cal A}$.
\end{lemma}

\begin{proof}
Let  $f$ be in ${\cal A}$, $x$ be a point in $X$ and $h$ be in the ideal ${\cal
I}(X)$.
Let $\gamma (t)$ be the integral curve of the Hamiltonian vector field
of $f$ with $\gamma (0) =x$.
Then $\gamma (t)$ is in $X$ and so $h(\gamma (t))= 0$ for all $t$.
Differentiation with respect to $t$ yields
$$
0 = \left. \frac{d}{dt}  \right|_0 h(\gamma (t)) = \{f, h\} (x)
$$
Thus $ \{f, h\}|_X =0$, i.e., $ \{f, h\}$ is in the ideal ${\cal I}(X)$.
\end{proof}
The Hamiltonian flows of invariant functions on $M$ preserve the fibers of the
momentum map $F$ (Noether's theorem). Therefore, by Lemma~\ref{lemma1} with
${\cal A} =
C^\infty (M)^G$ and $X = F^{-1}({\cal O})$ we have that
${\cal I}(F^{-1}({\cal O}))$ is a
Poisson ideal.  This proves that the smooth functions $C^\infty (M_{\cal O})$
on the reduced
space form a Poisson algebra.\\

\noindent
{\sc Remark}  More generally, we can define a sheaf of Poisson algebras
on the reduced space, a kind of structure sheaf.  An open set $U$ in the
 reduced space $M_{\cal O}$ is the quotient of
the intersection of the level set $F^{-1}({\cal O})$ with a $G$ invariant open
set $\tilde{U}$.  We define the Poisson algebra $C^\infty (U)$ by
$$
C^\infty (U) := \left. C^\infty (\tilde{U})^G \right| _{F^{-1}({\cal O})} .
$$

\noindent
{\sc Remark} Recall that if the action of the group $G$ on the manifold $M$ is
proper, then for a subgroup $H$ of $G$ the set of points $M_{(H)}$ of orbit
type $H$, i.e. the set of points with orbits isomorphic to $G/H$, is a
submanifold of $M$ (the definition of proper actions and some of their
properties are listed later).  Since the Hamiltonian flow of a $G$ invariant
function on
$M$ is $G$ equivariant, the manifolds $M_{(H)}$, $H<G$, are preserved by
the flows of invariant functions. Therefore for a subgroup $H$ of $G$
the ideal of invariant
functions  vanishing on the intersection  $F^{-1}({\cal O}) \cap
M_{(H)}$ is a Poisson ideal in the algebra of the invariant functions
and consequently defines a Poisson ideal in the algebra of smooth functions
$C^\infty (M_{\cal O})$ on the reduced space.\\

The Poisson bracket on the reduced space should allow us to write
down equations of motion for any $f \in C^\infty (M_{\cal O})$.  But
first we need to define what we mean by a smooth curve in a reduced
space.

\begin{definition}{\rm
A {\it smooth curve} $\gamma$ in a reduced space $M_{\cal O}$ is a continuous
map $\gamma :I \to M_{\cal O}$, $I$ an interval, such that for any smooth
function $h \in C^\infty (M_{\cal O})$ the function $h(\gamma (t))$ is a
smooth function on the interval $I$.

A {\it smooth flow} $\{\phi _s\}$ on $M_{\cal O}$ is defined similarly.
It is a one-parameter group of homeomorphisms $\phi_s: M_{\cal O} \to
M_{\cal O}$ such that for each $h\in C^\infty (M_{\cal O})$ and each
$s$, we have $h\circ \phi _s \in C^\infty (M_{\cal O})$, and
for each point $m\in M_{\cal O}$ the curve $s\mapsto \phi _s (m)$ is
a smooth curve.  }
\end{definition}
We are now in a position to define a Hamiltonian flow of a smooth function on a
reduced
space.

\begin{definition}{\rm
A {\em Hamiltonian flow of a smooth function $f$} on a reduced space
$M_{\cal O}$ is a smooth flow $\{\phi _s\}$ such that for any point
$m$ in $M_{\cal O}$ and any smooth function $h$ in $C^\infty(M_{\cal O})$
we have
\begin{equation}\label{eqHamilton}
	 \frac{d}{ds}h(\phi_s(m)) = \{f, h\} (\phi _s (m))
\end{equation}
where $\{\, , \,\}$ is the Poisson bracket on $C^\infty (M_{\cal O})$.
}
\end{definition}

This definition raises a problem. Since the reduced space $M_{\cal O}$
is not necessarily locally Euclidean, equation~(\ref{eqHamilton}) is
{\em not} in general a system of ordinary differential
equations in a coordinate-free notation.  Therefore the existence and
uniqueness of solutions of (\ref{eqHamilton}) needs to be addressed.

The existence is easy. The key fact is that the Hamiltonian flow of a
$G$ invariant function $\bar{f}$ on the original symplectic manifold
$M$ is $G$ equivariant.  Since the flow also preserves the level sets
$F^{-1}({\cal O})$ of the momentum map $F$, it descends to a flow on
the reduced space $M_{\cal O}$.  It is
now a formal exercise to check that this flow is smooth, and that it is a
Hamiltonian flow of the corresponding function $f\in C^\infty (M_{\cal
O})$ in the sense of the above definition (cf p.~389 in \cite{SL}).

The uniqueness is not to be expected without additional assumptions
about the topology of the reduced space.  Indeed, on a non-Hausdorff
manifold an integral curve of a vector field is not necessarily
unique.  One would expect non-uniqueness on any non-Hausdorff space.
The example of the irrational flow on the cotangent bundle of the two
torus considered above is quite instructive in this case.  Recall that
the reduced space at zero in the example is homeomorphic to the
product $({\Bbb T}^2/{\Bbb R}) \times {\Bbb R}$.  It is easy to see
that the smooth functions on this reduced space are simply the
functions that are constant on the first factor and smooth (in the
usual sense) on the second factor.  It follows that {\em any}
continuous flow on the product that fixes the points of the second
factor is smooth.  Since the induced Poisson bracket is zero, {\em
any} flow that fixes the points of the second factor is a Hamiltonian
flow of {\em any} smooth function. Thus a different set of ideas is needed
to make sense of non-Hausdorff reduced spaces.

\begin{lemma}
If the smooth functions on the reduced space separate points, then Hamiltonian
flows are unique.
\end{lemma}

\begin{proof}
Again we follow \cite{SL}, p.~389.  Suppose that $\phi_t$ and
${\psi_t}$ are two Hamiltonian flows on a reduced space $M_{\cal
O}$ generated by a function $f \in C^\infty (M_{\cal O}).$ Then, by
the chain rule, $\phi _{-t} $ is a flow of $-f$.  Since smooth functions
separate points, it is enough to show that for any function $h\in
C^\infty (M_{\cal O})$ and any point $m$ in the reduced space,
$$
	h(\psi _t (\phi _{-t}(m))) = h(m).
$$
However
$$
	\frac{d}{dt}h(\psi _t (\phi _{-t}(m))) = \{h, f\} (\psi _t (\phi _{-t}(m))) +
		\{h, -f\} (\psi _t (\phi _{-t}(m))) = 0.
$$
\end{proof}

At this point we make an assumption that will guarantee that functions
on the reduced space will separate points, namely that the action of
the symmetry group $G$ on the original manifold $M$ is {\em proper},
that is to say the map
$$
 G\times M \to M\times M, \quad (g, m)\mapsto
(g\cdot m, m)
$$
 is a proper map. Equivalently, an action of $G$ on
$M$ is proper if given two convergent sequences $\{m_n\}$ and
$\{g_n\cdot m_n\}$ in $M$ there exists a convergent subsequence
$\{g_{n_k}\}$ in $G$.

\subsection*{Digression:  properties of proper group actions}

We now list the properties of proper group actions that we will need in
the course of the paper.  The proofs of some properties are easy or are readily
available.  Other properties appear to be folklore and we will supply the
proofs in the appendix.

\noindent  1.
The isotropy group $G_m$ of any point $m$ in $M$ is compact; all orbits of $G$
in $M$ are closed and
embedded submanifolds.\\

\noindent  2. The orbit space $M/G$ is Hausdorff.\\

\noindent 3.  At every point $m \in M$ there exists a {\em slice } for the
action of $G$.  That is to say there is a ball $B$ about $0$ in the fiber
 $T_m M/ T_m G\cdot m$ of the normal bundle to the orbit through $m$
with $B$ invariant under the action of $G_m$ and an embedding $\phi : B \to
M$ with $\phi (0) = m$ such that the set $G\cdot \phi (B)$
is open in $M$ and the induced map
$$
	G\times _{G_m} B \to M, \quad [g, b] \mapsto g\cdot \phi (b)
$$
 is a diffeomorphism onto the image $G\cdot \phi (B)$.  Here $[g, b]$ denotes
the class of $(g, b) \in G\times B$ in the associated fiber bundle
$G\times _{G_m} B$ and $G_m$ is the isotropy group of $m$.
The $G_m$ invariant manifold $\phi (B)$ is a slice for the action of $G$
at $m$. \\

\noindent 4. There exists a $G$ invariant partition of unity subordinate to any
$G$ invariant open cover.
(We assume that the manifold $M$ is paracompact.)\\

\noindent 5. There exists on $M$ a $G$ invariant positive definite metric.\\

\noindent 6. Smooth $G$ invariant functions separate the orbits of $G$.\\

\noindent 7.  If $\omega $ is a $G$ invariant symplectic form on $M$ there
exists a $G$ invariant almost complex
structure $J$ {\em adapted} to $\omega$. That is to say, the bundle map $J$ is
symplectic, $\omega (J\cdot , J\cdot ) = \omega (\cdot , \cdot )$,
and the symmetric form $\omega (\cdot, J\cdot) $ is a positive definite metric.
A proof is provided in the appendix.\\

\noindent 8.  For any (compact) subgroup $H$ of $G$ the sets
$$
    M_{(H)} =
\{m \in M: G_m, \mbox{ the isotropy group of $m$, is conjugate to $H$ } \},
$$
$$
    M_{H} = \{m \in M: G_m \mbox{ is } H \},
$$
and
$$
    M^H = \{m \in M: G_m \mbox{ contains }H  \}= \mbox{ the set of points fixed
by $H$}
$$
are submanifolds of $M$ (this follows from the existence of slices ,
cf.\ fact~3.).  The manifold $M_{(H)}$ is called the manifold of points of
{\em orbit type } $(H)$.  Note also that the closure of $M_H$ is contained in
$M^H$ but need not equal all of $M^H$.\\

\noindent 9.  An equivariant version of the relative Darboux theorem holds:
\begin{theorem}[Relative Darboux]\label{thm.erDt}
 {\it Let $X$ be a submanifold of a manifold $Y$.
Let ${\omega }_0$ and ${\omega }_1$ be two symplectic forms on
$Y$ such that ${\omega }_0(x) = {\omega }_1(x)$  for each $x \in X $.
Then there exist neighborhoods ${\cal U}_0$ and ${\cal U}_1$ of $X$ and a
diffeomorphism $\psi :{\cal U}_0 \longrightarrow {\cal U}_1$ such that the
pull back of ${\omega }_1$ by $\psi $ is ${\omega }_0$ and
$\psi $ is the identity on $X$.

If a Lie group $G$ acts properly on
$Y$, preserves $X$, ${\omega }_0$, and ${\omega }_1$, then we can
arrange that the neighborhoods ${\cal U}_0$ and ${\cal U}_1$ are $G$-invariant
and that the diffeomorphism $\psi $ is $G$-equivariant.}
\end{theorem}
A proof is given in the appendix.\\

\noindent 10.  It follows from fact~9 that if $M$ is symplectic then the
manifolds $M_H$ and $M^H$ are symplectic as well.  The manifold $M_{(H)}$ is
usually not symplectic.

\section{Geometry of the reduced space }
The main goal of this section is to establish the following theorem.

\begin{theorem}\label{thm1}
Let $G$ be a Lie group acting properly and in a Hamiltonian way on
a symplectic manifold $(M, \omega )$ with a corresponding equivariant
momentum map $F : M \to {\frak g}^*$ and let ${\cal O}$ be
a locally closed coadjoint  orbit of $G$.  Then\\[3pt]
1. The reduced space $M_{\cal O} := F^{-1}({\cal O})/G$ is a locally finite
union of symplectic manifolds.  We will call these manifolds {\em symplectic
pieces}. \\[3pt]
2. The Hamiltonian flows of smooth functions preserve the decomposition of
the reduced space $M_{\cal O}$ into symplectic pieces.\\[3pt]
3.  The embedding of a symplectic piece into the reduced space $M_{\cal O}$ is
a Poisson map.
\end{theorem}

Observe that the condition of a coadjoint orbit being locally closed
is automatic for reductive groups and for their semidirect products
with vector spaces.  There is an example of a solvable group due to
Mautner (\cite{Puk}, p.~512) with non-locally closed coadjoint orbits
%
so the condition we are imposing is not vacuous.
Note also that the condition of the coadjoint orbit being locally closed
is precisely the condition that is necessary in order for the shifting trick to
make sense.  Since we want to read off the structure of the reduced space
from the corresponding slice representation on the original manifold we
will not use the shifting trick.

Theorem~\ref{thm1} has an important corollary.
\begin{corollary}
Suppose $M_{\cal O}$ and $N_{\cal O'}$ are two reduced spaces and
$\phi:  M_{\cal O} \to N_{\cal O'}$ a homeomorphism.  If the induced pull-back
map $\phi ^* C^\infty (N_{\cal O}) \to C^\infty (M_{\cal O}) $ is a Poisson
isomorphism then $\phi $ maps symplectic pieces to symplectic pieces.
\end{corollary}

\begin{proof}
On a connected
symplectic manifold the group generated by the time one Hamiltonian
flows of smooth functions acts transitively. It follows from this and
from assertion 2 of the theorem that connected components of the
symplectic pieces of a reduced space $M_{\cal O}$ are equivalence
classes of the relation: $x$ is equivalent to $y$ if and only if there
are smooth functions $f_1, \ldots, f_n \in C^\infty (M_{\cal O})$ such
that a composition of their time one flows maps $x$ to $y$.  Thus the
decomposition of a reduced space into symplectic manifolds is encoded
in the Poisson algebra $C^\infty (M_{\cal O})$ of smooth functions on
the reduced space.
\end{proof}

The corollary also allows us to define local isomorphisms of reduced spaces.
We will see in Theorem~\ref{thm.local}
 that all reduced spaces (under the two hypotheses above)
are locally isomorphic to a symplectic vector space reduced at zero with
respect to a linear action of
a compact group.  This in turn permits us to define abstract
``symplectic stratified spaces.''

Here is the strategy of the proof of Theorem~\ref{thm1}.  We will define the
terms and
 provide complete statements shortly.  Fix a point $x$
 in the preimage $F^{-1}({\cal O})$ of the coadjoint orbit.  Since the
 symplectic form $\omega $ on $M$ is $G$ invariant the $G$ orbit of
 $x$ is a submanifold of constant rank.  It follows by the constant
 rank embedding theorem of Marle \cite{marle:rank}, \cite{marle:model}
(see Theorem~\ref{thm.cret} below)
that a $G$ invariant neighborhood of the orbit
 $G\cdot x$ is symplectically determined by the restriction of the
 symplectic form to the orbit and by the symplectic normal bundle of
 the embedding $G\cdot x \hookrightarrow M$.  So if we can find a
 constant rank embedding of the orbit into some ``simple" Hamiltonian
 $G$ manifold $Y$, a neighborhood of the orbit in $Y$ is going to be
 symplectically isomorphic to a neighborhood of the orbit in $M$.  The
 manifold $Y$ is a kind of ``Darboux coordinates'' that take the
 action of $G$ into account.  We will then carry out our computations
 of the reduced space in $Y$.

\subsection*{Digression: constant rank embeddings}

Let $X$ be a submanifold of a symplectic manifold $(P, \tau)$.  For a point $x$
in $X$ the
{\em symplectic perpendicular} to the tangent space of $X$  at $x$ with respect
to the form
$\tau$ is the subspace
$$
T_xX ^\tau := \{v\in T_xP : \tau (x)(v, w) =0 \mbox{ for all } w\in T_x X\}.
$$
Together these subspaces define a subbundle $TX^\tau$ of the restriction $T_X
P$ of the tangent
bundle of $P$ to $X$.  Under the isomorphism
$$
T_X P \to T_X ^* P\quad (x, v) \mapsto (x, \tau (x) (v, \cdot))
$$
the bundle $TX^\tau$ is identified with the annihilator $TX^\circ $ of $TX$ in
$T_X^* P$.
Since
$TX^\circ = (T_X P/TX)^*$, the symplectic perpendicular $TX^\tau$ is
isomorphic, as an
abstract real vector bundle, to the normal bundle of $X$ in $P$.  In general
the form
$\tau (x)$ may be degenerate on $T_xX ^\tau$.  In fact, the quotient
$T_xX^\tau /(T_xX^\tau\cap T_xX)$ is isomorphic to a maximal symplectic
subspace of
$(T_xX^\tau, \tau(x))$.  If the dimension of this quotient is constant, i.e.,
if the distribution $TX\cap TX^\tau$ is a vector bundle, we say that the
embedding
$X\hookrightarrow (P, \tau)$ is of {\em constant rank}.  In this case the
quotient bundle
$$
N(X):= TX ^\tau/TX^\tau \cap TX
$$
is a symplectic vector bundle, the {\em symplectic normal bundle } of the
embedding
$X\hookrightarrow (P, \tau)$.  Note that $TX \cap TX^\tau $ is simply the
kernel of the
restriction of $\tau $ to $X$ and that as abstract vector bundles
$$
T_X P \simeq N(X)\oplus (TX^\tau \cap TX ) \oplus TX.
$$
So together the pull back $\tau |_X$ and the symplectic normal bundle $N(X)$
contain more information than the abstract normal bundle of
$X$ in $P$.  In fact the two pieces of data --- $\tau |_X$ and $N(X)$ --- {\em
uniquely
describe the symplectic form $\tau$ in a whole neighborhood of $X$}.  The
precise statement is
this.

\begin{theorem}[Uniqueness of constant rank embeddings]\label{thm.cret}
Let $(P, \tau)$ and $(P', \tau')$ be two symplectic manifolds.
Suppose $i:X \to (P,\tau)$ and $i':X \to (P',\tau')$ are two constant
rank embeddings with isomorphic symplectic normal bundles such that
$i^*\tau = (i')^*\tau '$.  Then there exist neighborhoods $U$ of
$i(X)$ in $P$ and $U'$ of $i'(X)$ in $P'$ and a diffeomorphism $\phi :
U \to U'$ such that $\phi \circ i =i'$ and $\phi ^* \tau ' = \tau$.

Furthermore, if $G$ is a Lie group that acts properly on $X$, $P$ and $P'$,
preserves the forms $\tau$ and $\tau'$ and if the embeddings $i$ and
$i'$ are $G$ equivariant, then $U$ and $U'$ can be chosen to be $G$
invariant and $\phi$ to be $G$ equivariant.
\end{theorem}

The proof of the theorem is given later in the appendix.
Since any two equivariant momentum maps differ by a constant vector, we also
have the
following corollary (we keep the notation of Theorem~\ref{thm.cret}).

\begin{corollary}\label{cor.m}
Suppose in addition that the actions of $G$ on $(P, \tau)$ and $(P', \tau')$
are
Hamiltonian with momentum maps $F: P \to {\frak g}^*$ and $F': P' \to {\frak
g}^*$.
If $F\circ i = F'\circ i'$ then $F'\circ \phi = F$.
\end{corollary}

\subsection*{}

This ends our digression and we now continue with the proof of
Theorem~\ref{thm1}.  Recall
that we have a Hamiltonian $G$ space $(M, \omega _M)$ with momentum map $F:M
\to {\frak g}^*$,
that $x$ is a point in $M$, $\alpha = F(x)$ and ${\cal O} =G\cdot \alpha$
is the coadjoint orbit through $\alpha $.  We want to
model a neighborhood of the orbit $G\cdot x$ in $M$ in order to understand
the structure of the quotient $F^{-1}({\cal O})/G$ near the orbit $G\cdot x$.

We have observed that $G\cdot x$ is a constant rank submanifold of $M$.  We
will now construct a symplectic manifold $(Y, \omega _Y)$ with a Hamiltonian
$G$ action, an equivariant momentum map $F_Y$ and an embedding
$i:G\cdot x \hookrightarrow Y$ such that $i^* \omega_Y=
\omega _M |_{G\cdot x}$, the symplectic normal bundle of $i$ is the same as of
$G\cdot x \hookrightarrow M$ and $F_Y (x) = \alpha $.  The constant
rank embedding theorem, Theorem~\ref{thm.cret}, would then guarantee
that there are neighborhoods $U$ of $G\cdot x$ in $M$, $U_Y$ of
$G\cdot x$ in $Y$ and a $G$ equivariant diffeomorphism $\phi : U_Y \to
U$ such that $\phi ^* \omega _M= \omega_Y$ and $\phi^* F = F_Y$.

\begin{proposition}\label{prop11}
Let $F: M \to {\frak g}^*$ be a momentum map for a Hamiltonian action of a Lie
group
$G$ on a symplectic manifold $(M,\omega_M)$. Then for any point $x\in M$ the
restriction of the ambient symplectic form $\omega_M$ to the orbit $G\cdot x$
equals the pullback by the momentum map $F$ of the symplectic form on the
coadjoint orbit through $F(x)$:
$$
\omega _M|_{G\cdot x} = F^* \omega _{G\cdot F(x)}|_{G\cdot x}
$$
where $\omega _{G\cdot F(x)}$ is the Kirillov-Kostant-Souriau (KKS) symplectic
form on the coadjoint
orbit of $F(x)$ in ${\frak g}^*$.
\end{proposition}

\begin{proof}
Let $\alpha = F(x)$.  For a vector $\xi \in {\frak g}$ let $\xi _M$ denote the
corresponding vector field on $M$ induced by the action of $G$
 and $\xi _{{\frak g}^*}$ the corresponding
vector field on ${\frak g}^*$ induced by the coadjoint action.
By the definition of the momentum map, we have for any $\xi, \eta \in {\frak
g}$
$$
\omega _M(x)(\xi _M (x), \eta _M (x)) = \langle \xi, dF_x (\eta _M (x))\rangle.
$$
By the equivariance of $F$, we have $dF_x (\eta _M (x)) = \eta_{{\frak g}^*}
(F(x))
=\eta_{{\frak g}^*} (\alpha )$.  Finally,
$$
\langle \xi, \eta_{{\frak g}^*} (\alpha )\rangle =\langle (-\ad \eta)\xi,
\alpha \rangle =
\langle [\xi, \eta ], \alpha \rangle = \omega _{G\cdot\alpha}(\xi_{{\frak
g}^*}(\alpha ),
\eta _{{\frak g}^*} (\alpha ))
$$
and we are done.
\end{proof}
{\sc Remark} The proposition allows us to take a more uniform view of
regular reduction and, in particular, to make sense of the shifting trick.
Suppose a momentum map $F: M \to \fg^*$ is transversal to a
locally closed orbit $G\cdot \alpha \subset \fg^*$ and that the action of
$G$ on the preimage $F^{-1}(G\cdot \alpha ) $ is free.  Then by the
proposition the form
$(\omega -F^* \omega _{G\cdot \alpha})|_{F^{-1}(G\cdot \alpha )}$
is basic and descends to a symplectic form on the quotient.  This form
on the quotient is the Marsden - Weinstein - Meyer reduced form.

\begin{corollary}\label{cor.k}
The symplectic perpendicular $T_x (G\cdot x)^{\omega _M}$ to the
tangent space at $x$ of  the orbit through $x$ intersects the tangent space
in the $G_\alpha$ directions:
$$
T_x (G\cdot x)^{\omega _M} \cap T_x (G\cdot x) =  T_x (G_\alpha \cdot x).
$$
Here $G_\alpha$ is the isotropy group of
$\alpha = F(x)$, the image of $x$ under the momentum map $F$.

\end{corollary}

Corollary~\ref{cor.k} says that the subspace
$({\frak g}_\alpha )_M (x) := \{\xi _M (x) :\xi \in {\frak g}_\alpha = Lie
(G_\alpha ) \} $
is the null space of the form $\omega _M(x)|_{T_x (G\cdot x)}$.
Note that since $F$ is equivariant, the isotropy group $G_x$ of $x$ is
contained in $G_\alpha$.

Choose a $G_x$ invariant splitting
\begin{equation}\label{eq1}
{\frak g} ={\frak g}_\alpha \oplus {\frak s}
\end{equation}
(we use the fact that $G_x$ is compact). Then $\omega _M(x)|_{\fs_M(x)} $ is
nondegenerate.  So $\fs_M(x)$ is a symplectic subspace
of $(T_x M, \omega _M(x))$ isomorphic to the tangent space of the
coadjoint orbit through $F(x)$.
(Note that $\omega _M (x)|_{\fs_M(x)}\simeq \omega _{G\cdot \alpha}
(\alpha )$, $\alpha = F(x)$.)\,  Also the symplectic perpendicular
$T_x (G\cdot x)^{\omega _M}$
contains $({\frak g}_\alpha )_M (x)$.  Pick a $G_x$ invariant
splitting $T_x (G\cdot x)^{\omega _M} =
({\frak g}_\alpha )_M (x) \oplus V$.  Then $V$ is a symplectic subspace
in   $(T_x M, \omega _M(x))$.
Note that $V$ is isomorphic to the quotient $T_x (G\cdot x)^{\omega _M}/
(T_x (G\cdot x)^{\omega _M}\cap  T_x (G\cdot x))$,
which is a typical fiber of the symplectic normal bundle of the
orbit $G\cdot x $ in $M$.  Let $\omega _V = \omega _M(x)|_V$.

Since $V$ and $\fs_M(x) $ are both symplectic  and $G_x$ invariant, the
symplectic perpendicular
$(V\oplus \fs_M(x))^{\omega _M}$ is also symplectic and $G_x$ invariant.
The symplectic perpendicular contains $({\frak g}_\alpha )_M (x)$,
which is null in it. Hence, by dimension count, this space is Lagrangian in the
symplectic perpendicular.
We conclude that $(T_x M, \omega _M (x))$ splits as a direct sum of three
symplectic  subspaces:
$$
 (T_x M , \omega _M(x)) = (V, \omega _V) \oplus ( \fs_M(x), \omega _{G\cdot
\alpha }(\alpha ))
		\oplus (({\frak g}_\alpha )_M (x)\oplus ( {\frak g}_\alpha )_M (x) ^*),
$$
and the splitting is $G_x$ invariant.
The symplectic form on the last summand is the canonical form on the product
of a vector space with its dual.

The tangent space at $x$ to the total space $Y$ of the associated bundle
$\pi: G\times _{G_x} ( (\fg_\alpha /\fg _x)^* \oplus V) \to G\cdot x$
(we think of the orbit $G\cdot x$ as being embedded
in the bundle as the zero section) is
$$
T_xY \simeq  T_x (G\cdot x) \oplus  (\fg_\alpha /\fg _x)^* \oplus V \simeq
	(\fs/\fg _x) \oplus (\fg _\alpha /\fg _x)\oplus (\fg_\alpha /\fg _x)^* \oplus
V
	\simeq T_x M .
$$
We now construct a closed $G$ invariant two form $\tau$ on the total space $Y$
of the associated bundle
such that
$$
 (T_x ( G\times _{G_x} [ (\fg_\alpha /\fg _x)^* \oplus V]), \tau (x) ) = (T_x
M, \omega (x)).
$$
The form $\tau $ is going to be the sum of three terms.
We construct the  first term $\tau _1$ by first pulling back by the momentum
map $F$
the KKS symplectic form
$\omega _{G\cdot \alpha }$ on the coadjoint orbit through $\alpha =F(x)$.  We
then pull it back by
the bundle projection map $\pi: Y \to G\cdot x$, so
 $\tau _1 = \pi ^* (F|_{G\cdot x})^* \omega _{G\cdot \alpha }$.
At the point $x$ the form $\tau _1$ is a nondegenerate two form on
the subspace $\fs/\fg _x\simeq \fs_M (x)$.

To construct the second and the third terms observe that the diagram
$$
\begin{array}{ccc}
Y= G\times _{G_x}  (\fg_\alpha /\fg _x)^* \oplus V& \longrightarrow &
{G\times _{G_x} V}\\
\downarrow & & \downarrow \\
 {G\times _{G_x}  (\fg_\alpha /\fg _x)^*}& \longrightarrow &     {G\cdot x}
  \end{array}
$$
commutes.  So we can think of $Y$ as a vector bundle over $G\times _{G_x} V$
with typical fiber $(\fg _\alpha /\fg _x)^*$ or as a vector bundle over
$G \times_{G\cdot x} (\fg _\alpha /\fg _x)^* $ with typical fiber $V$.
Therefore a form on the total space of $G\times _{G_x} V$ or of
$G\times_{G\cdot x} (\fg _\alpha /\fg _x)^* $ may be thought of as a form on
the
manifold $Y$.

To construct the second term $\tau _2$ we embed $G\times _{G_x} (\fg _\alpha /
\fg _x )^* $ into the cotangent bundle of the orbit $G\cdot x$.  The $G_x$
equivariant splitting $\fg = \fs\oplus \fg_\alpha $ chosen above gives rise
to a $G_x$ equivariant projection
$$
\fg \to \fg_\alpha ,
$$
which induces an embedding
$$
 j : (\fg _\alpha /\fg _x)^* \hookrightarrow (\fg / \fg _x)^*
$$
and thereby an embedding $j$ of the associated vector bundles
$$
 j: G\times _{G_x} (\fg _\alpha / \fg _x )^* \to
 G\times _{G_x} (\fg / \fg _x )^* \simeq T^* (G\cdot x).
$$

The pull-back by $j$ of the canonical symplectic form
$\omega _{T^* (G\cdot x)}$ on the cotangent bundle of the orbit $G\cdot x$
is  a closed two form on $ G\times _{G_x} (\fg _\alpha / \fg _x )^*$
hence gives rise to a closed two form $\tau _2$ on $Y$.

The construction of the third term $\tau _3$ is an example of minimal coupling
of Sternberg.  We first refine the splitting
 $\fg = \fs\oplus \fg_\alpha $
to a $G_x$ equivariant splitting
$$
\fg = \fg_x \oplus \fm \oplus \fs
$$
(with $\fg_\alpha  = \fg_x \oplus \fm $).  Let $A_0 : \fg \to \fg_x $ be the
corresponding $G_x$ equivariant projection.  It defines a left $G$ invariant
$\fg_x$-valued one form $A$ on $G$.  The form $A$ is a connection one form for
the
principal $G_x$ bundle $G\to G\cdot x$.  Let $F_V : V \to \fg _x ^* $ denote
the momentum map arising from the linear symplectic action of $G_x$ on
$(V, \omega _V)$.  Consider the following two form on $G\times V$:
$$
	\tilde{\tau }_3 = d\langle A, F_V \rangle + \omega _V.
$$
The form is $G_x$ invariant.  It is not hard to check that it is basic
for the projection $G\times V \to G\times _{G_x} V$.  Denote the corresponding
closed two form on $G\times _{G_x } V$ and hence on $Y$ by $\tau _3$.
Note that the value of $\tau _3 $ at $x$ on $V\subset T_xY =T_x M$
is the form $\omega_V$
 and the restriction of $\tau _3(x)$ to the other two
summands is zero. (This is because $F_V$ is a homogeneous quadratic map on $V$,
so $F_V (0)  = 0$ and $dF_V (0) =0$. Consequently
 at a point $(g, 0)\in G \times V$ we have
$
d\langle A, F_V \rangle (g,0) = \langle dA, F_V \rangle (g, 0) + \langle
A \wedge dF_V \rangle (g, 0) = 0 + 0 = 0.
$)

We conclude that $(\tau _1 + \tau _2 + \tau _3)(x) = \omega (x)$.
Let $\tau $ be the sum $ \tau _1 + \tau _2 + \tau _3 $.
By construction $\tau $
is a closed $G$ invariant two form on the total space of $Y= G\times _{G_x}
((\fg_\alpha /\fg _x)^* \oplus V)$ which is  non-degenerate at the points of
the zero section $G\cdot x$.  Thus $\tau $ is non-degenerate in some
($G$ invariant) neighborhood $Y_0$ of the zero section.  Note that by
 construction we have $\tau |_{G\cdot x} = F^* \omega _{G\cdot \alpha } =
\omega |_{G\cdot x}$ and the symplectic normal bundle of the embedding
$G\cdot x \hookrightarrow (Y_0, \tau )$ is $G\times _{G_x} V$ which is
the symplectic normal bundle of the embedding
$G\cdot x \hookrightarrow (M, \omega )$.
The constant rank embedding theorem says that if we shrink the
neighborhood $Y_0$ a bit more, we will have a $G$ equivariant map
$\psi: Y_0 \to M$ which is a diffeomorphism onto its image and has the
properties that $\psi |_{G\cdot x}$ is the identity map and $\psi ^* \omega =
\tau$.  That is to say,  a neighborhood of the zero section  $G\cdot x$ in
the associated bundle
$(G\times _{G_x} [ (\fg_\alpha /\fg _x)^* \oplus V], \tau)$ is the
promised ``Darboux coordinate patch adapted to the action of the Lie group
 $G$.''

Our next step  is compute a momentum map $F_Y$ for the action of $G$ on
$(Y, \tau)$.  The requirement that $F_Y (x) = F(x) = \alpha$ would then
ensure that $F_Y = F\circ \psi$, where $F:M \to \fg ^*$ is the original
momentum map.  This would finally allow us to get on with computing the
the reduced space $F^{-1} (G\cdot \alpha )/G$.

A momentum map is traditionally defined for actions that preserve {\em
nondegenerate} two forms.  One can extend this definition to actions that
preserve arbitrary two forms as long as the contractions of the vector fields
induced by the action with the form are exact.  For example, consider
the form $\tau _1 = \pi^* (F|_{G\cdot x})^* \omega _{G\cdot \alpha}$ on the
manifold $Y$, where $\pi : Y \to G\cdot x$ is the vector bundle projection.
For any $\xi \in \fg$ we have
$$
\xi _Y \lefthook\tau _1 =d \langle \xi, \pi^* (F|_{G\cdot x})\rangle .
$$
Hence $F_1 = \pi^* (F|_{G\cdot x})$ is a momentum map for the action of
$G$ on $(Y, \tau _1)$.  Note that $F_1 ([g, \lambda , v]) = g\cdot \alpha$
where $[g, \lambda , v]$ is the class of $(g, \lambda , v) \in G\times
(\fg _\alpha /\fg _x)^* \times V$ in the associated bundle $Y$.

Similarly since the lifted action  of $G$ on the cotangent bundle
$T^* (G\cdot x)= G\times _{G_x} \fg _x^\circ$ ($\fg _x^\circ$ is the
annihilator of $\fg _x$ in $\fg^*$) is Hamiltonian with momentum
map sending the class $[g, \lambda ]$ of
 $(g, \lambda )\in G\times \fg_x ^\circ$ to $g\cdot \lambda $, the map
$F_2 : G\times _{G_x} [ (\fg_\alpha /\fg _x)^* \oplus V] \to \fg ^*$
sending the class $[g, \lambda , v]$ to $g\cdot j(\lambda )$ is a momentum
map for the action of $G$ on $(Y, \tau _2 )$.  Recall that
$j: (\fg _\alpha /\fg_x)^* \to (\fg/\fg_x)^*$ is defined by the choice of
the splitting $\fg = \fg_\alpha \oplus \fs$ made earlier.

Let us also compute a momentum map for the action of $G$ on $(Y,\tau_3)$.
Note first that the action of $G$ on $G\times V$ given by
$g\cdot ( a, v) = (ga, v)$ preserves the form $\tilde {\tau }_3 = d\langle
A, F_V\rangle + \omega _V$.   So for $\xi \in \fg$ the induced vector field
$\xi_{G\times V} =\xi _G$ is a right invariant vector field on $G$ and
$$
0 = (\xi _G\lefthook d + d \xi _G\lefthook )\langle A, F_V \rangle \\
= \xi _G\lefthook d \langle A, F_V \rangle
	+ d \langle A (\xi _G), F_V \rangle .
$$
Since $\xi _G $ is right invariant and $A$ is left invariant, $A(\xi _G)(g) =
A_0 (Ad(g^{-1})\xi)$.  It follows that a momentum map for the action of
$G$ on $(G\times V, \tilde {\tau }_3)$ is given by
$(g, v) \mapsto g\cdot i (F_V (v))$ where $i:\fg _x ^* \to \fg^* $ is the
transpose of the projection $\fg \to \fg_x$. Therefore
$
	F_3 : (Y, \tau_3 ) \to  \fg  ^*,
$
a momentum map for the action of $G$  on $(Y, \tau_3 )$, is given by
$$
	F_3 ([g, \lambda , v])= g\cdot i(F_V (v)).
$$

The upshot of these computations is that $F_Y = F_1 + F_2 + F_3$ is a momentum
map for the action of $G$ on $(Y, \tau = \tau_1+ \tau _2 + \tau_3)$, that is to
say
$$
	F_Y([g, \lambda , v])= g\cdot ( \alpha + j(\lambda ) + i(F_V (v)))
$$
where  $i: \fg ^*_x \hookrightarrow \fg  ^*$ and
$j: (\fg _\alpha /\fg _x)^* \hookrightarrow \fg _x ^\circ $ are induced by an
$G_x$ equivariant splitting
$$
\fg = \fg _\alpha \oplus \fs= \fg _x \oplus \fm \oplus \fs.
$$
\vspace{13pt}

The proposition below is a key computation.
\begin{proposition}\label{prop13}
Assume that the coadjoint orbit through $\alpha = F(x)$ is locally closed.
Then for a small enough neighborhood $Y_0$ of the orbit $G\cdot x$ in the
model space $Y$, the intersection of the set  $F_Y^{-1}(G\cdot \alpha)$ with
the neighborhood $Y_0$ is of the form
$$
F_Y^{-1}(G\cdot \alpha) \cap Y_0 = \{[g, \lambda , v] \in Y_0:
\lambda =0 \mbox{ and } F_V (v) = 0 \}.
$$
\end{proposition}
\begin{proof}
We write $$F_Y : G\times _{G_x} ((\fg_\alpha /\fg _x)^* \times V) \to \fg
^*,\quad
[g, \lambda ,v]\mapsto g\cdot (\alpha + j(\lambda ) + i(F_V(v)))$$ as a
composition of two maps:
$$
b: G\times _{G_x} ((\fg_\alpha /\fg _x)^* \times V) \to
	 G\times _{G_x} \fg _\alpha ^*, \quad
 {[g, \lambda ,v]} \mapsto [g, \lambda + i_{\alpha} (F_V(v))]
$$
and
$$
{\cal E}: G\times_{G_x} \fg _\alpha ^* \to \fg^*, \quad
	{[g, \nu ]} \mapsto g\cdot (\alpha + k(\nu)).
$$
Here $i_{\alpha}: \fg _x ^* \hookrightarrow \fg_\alpha ^*$ is the $G_x$
equivariant embedding defined
by the $G_x$ equivariant splitting $\fg_\alpha = \fg_x \oplus \fm$ chosen
earlier,
$(\fg_\alpha /\fg_x)^* $ is identified with the annihilator of $\fg _x$ in $\fg
_\alpha ^*$ and
$k:\fg_\alpha ^* \hookrightarrow \fg^*$ is the $G_x$ equivariant embedding
corresponding to
the splitting $\fg = \fg_\alpha \oplus \fs$.  At the points of the form $[g,0]$
the map ${\cal E}$
is a submersion.  By assumption the coadjoint orbit $G\cdot \alpha $ is
embedded.  Therefore
the preimage ${\cal E}^{-1} (G\cdot \alpha )$ is an embedded submanifold of
$G\times _{G_x} \fg^* _\alpha $ of codimension $\dim \fg_\alpha$.  It follows
that the zero section
of $G\times _{G_x} \fg^* _\alpha $ is a collection of
connected components of the preimage of the orbit.
Since the preimage is embedded, there is a neighborhood ${\cal U}$ of the zero
section such that
${\cal E}^{-1} (G\cdot \alpha )\cap {\cal U}$ is the zero section.  Let $Y_0 =
b^{-1}({\cal U})$.
Then
$$
F_Y ^{-1} (G\cdot \alpha )\cap Y_0 = b^{-1} ({\cal U} \cap {\cal E}^{-1}
(G\cdot \alpha )) =
	Y_0 \cap b^{-1} (\mbox{zero section}).
$$
Clearly, $b^{-1} (\mbox{zero section}) =\{[g,0,v]: F_V (v) =0\}$ and we are
done.
\end{proof}

\begin{corollary}\label{cor14}
For any subgroup $H <G$ and any $\alpha \in \fg^*$ with the orbit $G\cdot
\alpha $ locally
closed, the set
$F^{-1}(G\cdot \alpha ) \cap M_{(H)}$ is a submanifold of $M$ of constant rank,
the quotient
$(M_\alpha )_{(H)}:= (F^{-1}(G\cdot \alpha ) \cap M_{(H)})/G$ is a symplectic
manifold and
the inclusion $(M_\alpha )_{(H)}\hookrightarrow M_\alpha := F^{-1}(G\cdot
\alpha )/G$ is a Poisson map.
\end{corollary}

\begin{proof} Let $x$ be a point in the intersection $F^{-1}(G\cdot \alpha )
\cap M_{(H)}$. It is no loss of generality to assume that the isotropy group
$G_x $ of $x$ is  $H$.
Recall that there exists a $G$ invariant neighborhood $U$ of the orbit
$G\cdot x$ in $M$, a $G$ invariant neighborhood $Y_0$ of the zero section in
$Y= G\times _{G_x} ((\fg_\alpha /\fg_x )^*
\times V)$ and a $G$ equivariant diffeomorphism $\psi : Y_0 \to U$  such that
$\psi ^* \omega = \tau |_{Y_0}$ where $\tau $ is the closed
two form on $Y$ constructed
earlier.  It follows that $\psi$ descends to a homeomorphism
$\psi _\alpha :  (Y_0)_\alpha \to U_\alpha$.  Here $(Y_0)_\alpha :=
(F_Y^{-1} (G\cdot \alpha ) \cap Y_0 ) /G$ and $U_\alpha := (U\cap F^{-1}(G\cdot
\alpha ))/G$
are the reduced spaces.  Moreover by construction the pull-back map
$\psi_\alpha ^*: C^\infty ((Y_0)_\alpha ) \to C^\infty (U_\alpha )$ is an
isomorphism of Poisson algebras, where
$C^\infty (U_\alpha ):= C^\infty (U)^G |_{F^{-1}(G\cdot \alpha )}$ etc.
Therefore it is enough to prove the statements of the corollary for the action
 of $G$ on $(Y_0, \tau)$.

It is convenient to ignore the distinction between the total space $Y$ and the
neighborhood $Y_0$ of the zero section in $Y$. We note for future reference
that the embedding of the symplectic
slice $V$ into the model space $Y$ given by $v\mapsto [e,0,v]$, where $e$ is
the identity element of $G$, is symplectic, i.e., $\tau |_V = \omega _V$.
It is not hard to see that
$$
Y_{(H)} =  G\times _H [(\fg _\alpha /\fh)^* \times V]^H
$$
where $[(\fg _\alpha /\fh)^* \times V]^H$ denotes the subspace of $H$ fixed
vectors (remember that $G_x = H$).
It follows that $F_Y^{-1}(G\cdot \alpha ) \cap Y_{(H)} = G\times _H V^H\simeq
G/H \times V^H$.
Since $\tau |_{G/H \times V^H} = \omega |_{G\cdot x} + \omega _V | _{V^H}$ and
 $V^H$ is a symplectic
subspace of $V$, we conclude that $F_Y^{-1}(G\cdot \alpha ) \cap Y_{(H)}$ is a
submanifold
of $(Y, \tau)$ of constant rank and that the quotient $(F_Y^{-1}(G\cdot \alpha
) \cap Y_{(H)})/G$
is diffeomorphic to $V^H$.

Therefore $F^{-1}(G\cdot \alpha ) \cap M_{(H)}$
is a submanifold of $(M, \omega )$ of constant rank and the quotient
$(F^{-1}(G\cdot \alpha ) \cap M_{(H)})/G$ is a manifold.    We also see that
the form
$(\omega - F^* \omega _{G\cdot \alpha })|_{F^{-1}(G\cdot \alpha ) \cap
M_{(H)}}$ is basic and
descends to a form on the base locally isomorphic to $\omega _V |_{V^H}$.
Thus $(M_\alpha )_{(H)}$ is a symplectic manifold.

Finally observe that if a function $f$ on $Y$ is $G$ invariant then
its restriction to
${F_Y^{-1}(G\cdot \alpha ) \cap Y_{(H)}} = {G/H \times V^H} $ is completely
determined by  by its restriction to ${V^H}$. It follows that the map
$(M_\alpha )_{(H)}\hookrightarrow M_\alpha $ is Poisson.
\end{proof}

To finish the proof of Theorem~\ref{thm1} we need to show that the
decomposition of the reduced
space $M_{G\cdot \alpha } = M_\alpha $ into symplectic pieces is locally finite
and that the
Hamiltonian flows of smooth function on the reduced space preserve the
decomposition.
The local finiteness of the decomposition of the reduced space follows from the
local finiteness
of the decomposition of the original manifold $M$ into orbit types which in
turn follows
from the existence of slices for proper group actions.  The fact that the
Hamiltonian flows
preserve the symplectic pieces of the reduced space is a consequence of the
fact that Hamiltonian
flows of $G$ invariant functions on $M$ are $G$ equivariant and hence preserve
the orbit types.
This concludes the proof. \hspace*{\fill}$\Box$\\

We now refine Theorem~\ref{thm1} in two different ways.
The first refinement is a theorem that
describes more precisely how the symplectic pieces fit
together.
The second one is a theorem that shows that the symplectic pieces can also be
obtained by regular reduction, thus providing a way to reconstruct the
reduced dynamics.

\begin{theorem}\label{thm.local}
Let $G$ be a Lie group acting properly and in a Hamiltonian fashion on
a symplectic manifold $(M, \omega )$ with a corresponding equivariant
momentum map $F : M \to {\frak g}^*$.  Let $x$ be a point in the manifold
$M$ and $\alpha = F (x)$.  Assume that the coadjoint orbit
$G\cdot \alpha$ is locally closed. Then locally the reduced space
$M_{\alpha} = F^{-1} (G\cdot \alpha )/G$ is isomorphic to a
neighborhood of the image of the origin in the space obtained by
reduction at $0$ of the symplectic slice
$V := T_x( G\cdot x)^\omega /(T_x (G\cdot x)^\omega \cap T_x (G\cdot x)) $
 through the point $x$ with respect to the isotropy group $G_x$ of $x$.
\end{theorem}

\begin{proof}
As before it is enough to prove the theorem at the point $x$ for the model
space $(Y, \tau)$, where
$Y= G\times _{G_x} ((\fg_\alpha /\fg _x)^* \times V)$.
  Of course, as was mentioned previously, the form $\tau$
is only nondegenerate on some neighborhood of the zero section $G\cdot x$, so
the Poisson bracket is only defined in that neighborhood.  It would be a
notational nightmare to keep track of the neighborhood so we will again
pretend that $\tau $ is nondegenerate everywhere on $Y$.

Recall that the embedding of the symplectic
slice $V$ into the model space $Y$ given by $v\mapsto [e,0,v]$,
where $e$ is the identity element
of $G$, is symplectic, i.e., $\tau |_V = \omega _V$.
We now prove that the restriction $C^\infty (Y) \to C^\infty (V)$ induces
an isomorphism of Poisson algebras $C^\infty (Y_\alpha) \to C^\infty (V_0)$
where $C^\infty (Y_\alpha):= C^\infty (Y)^G|_{F_Y ^{-1}(G\cdot \alpha) }$
and $C^\infty (V_0):= C^\infty (V)^{G_x}|_{F_V^{-1}(0)}$ are the algebras
of smooth functions on the corresponding reduced spaces.
It is easy to see that restriction to $V$ defines a surjective map
$C^\infty (Y)^G \to C^\infty (V)^{G_x}$, $f\mapsto f|_V$.
  Indeed, any $G$ invariant function
on $ G\times _{G_x} ((\fg_\alpha /\fg _x)^* \times V)$ restricts to an
$G_x$ invariant function on $(\fg_\alpha /\fg _x)^* \times V$ hence to a
$G_x $ invariant function on $V$.  Conversely, any $G_x$ invariant function
on $V$ extends trivially to a $G\times G_x$ invariant function on
$ G\times (\fg_\alpha /\fg _x)^* \times V$, so the restriction map
$C^\infty (Y)^G \to C^\infty (V)^{G_x}$ is surjective.  In fact the same
argument shows that the map
$C^\infty (G\times_{G_x}(\{0\} \times V))^G \to C^\infty (V)^{G_x}$ is
bijective.  Since $F_Y ^{-1}(G\cdot \alpha) \cap V = F_V ^{-1} (0)$, it follows
that the map
$$
 C^\infty (Y)^G|_{F_Y ^{-1}(G\cdot \alpha) } \to
		C^\infty (V)^{G_x}|_{F_V^{-1}(0)}
$$
induced by restriction to $V$ is bijective as well.
The bijection is a Poisson map because  $\tau |_V = \omega _V$.
\end{proof}

Theorem~\ref{thm.local} shows that the
decomposition of a reduced space into symplectic pieces is well behaved.
The reason for this good behavior is that the decomposition of a vector space
reduced at zero with respect to a linear action of a compact group
forms a Whitney regular stratification.  We now present the details.

We recall the discussion in \cite{SL}.  Let $(V, \omega _V)$ be
a symplectic vector space, $K$ a compact Lie group  and $K \to {\sl Sp}\, (V,
\omega _V)$ a symplectic representation of $K$.  As was mentioned earlier
the $K$ momentum map $F_V: V \to {\frak k}^*$ that sends zero to zero
is a quadratic polynomial.  The reduced space at zero $V_0 = F_V^{-1} (0)/K$
can be described as a semi-algebraic set.  To this end consider the algebra
${\Bbb R}[V]^K$ of $K$ invariant polynomials on $V$.  It is well known that
the algebra is finitely generated.  A deep result due to G. Schwarz \cite{Sch}
asserts that the algebra of invariant functions $C^\infty (V)^K$  is also
finitely generated in the following sense.
Let $p_1, \ldots, p_n$ be a minimal set of generators of the algebra of
invariant polynomials and let $p:V \to {\Bbb R}^n$ be given by
$p(v) = (p_1(v), \ldots, p_n (v))$.  Schwarz's theorem asserts
that the smooth invariant functions on $V$ are the compositions of smooth
functions on ${\Bbb R}^n$ with the invariant map $p$, i.e.,
$$
C^\infty (V)^K = p^* C^\infty ({\Bbb R}^n).
$$
Since $K$ invariant functions separate orbits the induced map $\bar{p} :V/K \to
{\Bbb R}^n$ is injective.  In fact it is a proper embedding (see for example
\cite{Bie}) and the image $\bar{p}(V/K) = p(V)$ is, by the Tarski - Seidenberg
theorem, a semi-algebraic subset of ${\Bbb R}^n$.

It is also easy to see that the map $\bar{p}$ embeds the reduced space $V_0$
as a semi-algebraic subset.  Indeed, let $||\cdot ||$ be a norm on
 the dual of the Lie algebra ${\frak k}^*$ defined by  a $K$ invariant
inner product.  Then $||F_V ||^2$ is an invariant polynomial on $V$.  So there
is a polynomial $f$ on ${\Bbb R}^n$ such that $||F_V ||^2 = f\circ p$.
Since $(||F_V ||^2)^{-1} (0) = F_V ^{-1} (0)$ we have
$\bar{p}(V_0) = \{f = 0\} \cap \bar{p}(V)$.  Thus $ \bar{p}(V)$ is a
semi-algebraic set.  Note that in complete analogy with Schwarz's theorem
the embedding map $\bar{p}:V_0 \to {\Bbb R}^n$ induces a surjective map
$\bar{p}^* : C^\infty ({\Bbb R}^n)\to C ^\infty (V_0 )$,
where as before $C^\infty (V_0 )$ denotes the algebra of smooth functions on
the reduced space, $C ^\infty (V_0 ) = C^\infty (V)^K |_{F_V^{-1} (0)}$.

It was shown in \cite{SL} that $\bar{p}$ embeds symplectic pieces of the
reduced space $V_0$ into smooth submanifolds of ${\Bbb R}^n$.  This defines
a decomposition of the semi-algebraic set $\bar{p}(V_0)$ into smooth
manifolds.  It was also shown (loc.\ cit.) that this decomposition of
$\bar{p}(V_0)$ satisfies the Whitney regularity condition.

Thus Theorem~\ref{thm.local} asserts that the decomposition of reduced spaces
into symplectic pieces defined by orbit type is a stratification
(in a technical sense) and that the stratification is locally Whitney regular.
This excludes more pathological singular spaces such as a cone over the
integers (not locally finite) or
 the set
$$
X = \{(x,y) \in {\Bbb R}^2:x=0,\,\, -1\leq y \leq 1 \}
\cup \{(x,y) \in {\Bbb R}^2: x>0,\, y =\sin \frac{1}{x}\},
$$
which is connected but not path connected.

The next theorem is useful in reconstructing the original dynamics from the
dynamics in the reduced space.

\begin{theorem}\label{thm.lift}
 Let $G$ be a Lie group acting properly and in a Hamiltonian
way on a symplectic manifold $(M, \omega )$ with a corresponding
equivariant
momentum map $F : M \to {\frak g}^*$.  Let $x$ be a point in the manifold
$M$, $H$ the isotropy group of $x$  and $\alpha = F (x)$.

The manifold
$$
	M_H:= \{m\in M: G_m = H \}
$$
is a symplectic submanifold of $M$. The normalizer $N$ of $H$ in $G$ preserves
$M_H$ and the quotient group  $L= N/H$ acts freely. This action of $L$ on $M_H$
is Hamiltonian and the reduced space $(M_H)_{\alpha _0}$,
for an appropriate vector $\alpha _0\in \fl^*$ is isomorphic to the symplectic
piece
$(M_\alpha )_{(H)}= (F^{-1}(G\cdot \alpha ) \cap M_{(H)} )/G$
provided the orbit $G\cdot \alpha$ is locally closed.
\end{theorem}
{\sc Remark}  Note that the action of $L$ on $M_H$ is free by construction,
so the manifold $(M_H)_{\alpha _0}$ is obtained by regular Marsden - Weinstein
- Meyer reduction.

\begin{proof}
We show first that the action of $L$ on $(M_H, \omega |_{M_H})$ is Hamiltonian.
 It is no loss of generality to assume that the manifold $M_H$ is connected.
Note that by the definition of $N$ the Lie algebra $\fn$ satisfies
$$
	\fn = \{X\in \fg: [X, \fh ] \in \fh\},
$$
where $\fh = \Lie(H)$.  Therefore the image of $\fl = \Lie (L)\simeq \fn /\fh$
under the embedding $\fn/\fh \hookrightarrow \fg/\fh$ is the set of $H_0$
fixed vectors $(\fg/\fh)^{H_0}$,  where $H_0$ is the identity component
of $H$.  It follows that the dual projection $(\fg/\fh)^* \to (\fn/\fh)^*$
restricts to an {\em isomorphism}
$\pi _{\fl} : [(\fg/\fh)^*]^{H_0} \to (\fn/\fh)^* = \fl^*$.  Therefore $\fl^* $
is
naturally isomorphic to $(\fh^\circ )^{H_0}$, the $H_0$ fixed vectors in
the annihilator $\fh^\circ$ of $\fh$ in $\fg^*$.

Recall that $x$ is a point in $M_H$ and  $\alpha$ is its image under the
$G$ momentum map $F$.  We claim that the image of $M_H$ under $F$ lies in
the affine plane $(\fh^\circ)^H + \alpha $.  Indeed, since $M_H$ is
pointwise fixed by $H$ and $F$ is equivariant, the image $F(M_H)$ is also
pointwise fixed by $H$.  Also, for any vector $\xi \in \fh$, any point
$y\in M_H$ and any tangent vector $v\in T_y M_H$ we have
$$
\langle \xi, dF_y (v)\rangle = \omega (y) (\xi _M (y) , v) = 0
$$
since $\xi _M(y) = 0$ for all $y\in M_H$.  Thus $dF_y (T_y M_H) \subset
\fh^\circ $ for all $y\in M_H$ and so $F(M_H) \subset \fh ^\circ + \alpha$
since $F(x) = \alpha $.

We conclude that the map $F_L := \pi _{\fl} \circ (F|_{M_H})$ is a
momentum map for the action of $L$ on $(M_H , \omega |_{M_H})$.

Since $H$ is closed in $G$, its normalizer $N$ is also closed in $G$,
so the action of $N$ on $M_H$ (and hence of $L$) is proper.  Therefore the
reduced space $(M_H)_{\alpha _0} := F_L ^{-1} (L\cdot \alpha _0)/L$ is
a symplectic manifold.  As was mentioned before, Proposition~\ref{prop11}
allows us the following description of the reduced symplectic structure on
$(M_H)_{\alpha _0}$.  The form $\omega _{M_H} :=
\omega |_{M_H}$ is not basic when restricted to the principal $L$ bundle
$F_L^{-1}(L\cdot \alpha_0) $, but the difference
$(\omega _{M_H} - F_L^* \omega _{L\cdot \alpha_0})
	|_{F_L^{-1}(L\cdot \alpha_0)}$
is basic and descends to the reduced symplectic form on $(M_H)_{\alpha _0}$
(here $\omega _{L\cdot \alpha_0}$ is the KKS symplectic form
on the coadjoint orbit $L\cdot \alpha_0 \subset \fl^*$).

We are now ready to prove the main claim of the theorem: that the manifold
$(M_\alpha ) _{(H)} := (F^{-1}(G\cdot \alpha )\cap M_{(H)})/G$ is
symplectically diffeomorphic to $(M_H)_{\alpha _0}$.  Note first that
$F_L^{-1}(L\cdot \alpha _0) = F^{-1} (N\cdot \alpha )\cap M_H$.  So to
establish the diffeomorphism, it is enough to show that
$$
(F^{-1} (N\cdot \alpha )\cap M_H)/N \simeq
	(F^{-1} (G\cdot \alpha )\cap M_{(H)})/G.
$$
We computed the right hand side locally in the proof of Corollary~\ref{cor14}.
We now compute the left hand side locally using the same model $(Y_0, \tau)$.
(As before we will ignore the distinction between the neighborhood $Y_0$ and
the whole space $Y$.)\,
We will see that locally the reduced space $(M_H)_{\alpha _0}$ is modeled
by the vector space $(V^H, \omega _V|_{V^H})$.  This will prove the theorem.

The manifold $Y_H$ of points in $Y$ with isotropy group $H$ is equal to
$N\times _H [(\fg_\alpha /\fh )^* \times V]^H$.  An argument
similar to the proof of Proposition~\ref{prop13} (the factoring of the map
$F_Y$ through two maps) shows that
$$
Y_0\cap F_Y^{-1} (N\cdot \alpha ) = Y_0 \cap (N\times _H F^{-1}_V (0)).
$$
It follows that
$$
{(Y_0)}_H \cap  F_Y^{-1} (N\cdot \alpha ) = N\times _H V^H \simeq L \times V^H.
$$
Therefore, locally,
$$
(M_H)_{\alpha _0} \simeq (V^H, \omega _V|_{V^H})
$$
and we are done.
\end{proof}

\subsection*{Actions of compact Lie groups: coadjoint directions don't matter}

The proofs of theorems \ref{thm1},  \ref{thm.local} and \ref{thm.lift} were
based on Marle's constant rank embedding theorem \cite{marle:rank},
 \cite{marle:model}.  However for compact symmetry groups we can
also use a local normal form theorem due to Guillemin and Sternberg \cite{GS1}.
This normal form is based on the idea of symplectic cross-sections.  It allows
us to restrict our attention to the smallest symplectic submanifold containing
a given fiber of the momentum map.  This reduces the number of the degrees
of freedom and the dimension of the symmetry group.  As a result, for compact
symmetry groups we only need to deal with reduction at zero values of the
momentum maps which is described in \cite{SL}.

The symplectic cross-section theorem can be stated as follows.

\begin{theorem}[{\rm \cite{GS2}, Theorem~26.7}] \label{thm:cross-section}
Let $(M, \omega )$ be a Hamiltonian $G$ space with  momentum map
$F : M\to {\fg}^*$.  Let $S$ be a submanifold of ${\fg}^*$ passing
through a point $\alpha \in {\fg}^*$ satisfying $T_\alpha S \oplus
T_\alpha G\cdot\alpha = {\fg}^*$ and $S$ is $G_\alpha $ invariant.
Then for a small enough neighborhood $B$ of $\alpha $ in $S$ the
preimage $F^{-1} (B)$ is a {\em symplectic} submanifold of $M$.

Moreover if we choose  $B$ to be  $G_\alpha $ invariant,
then the action of $G_\alpha $ on
$F^{-1} (B)$ is Hamiltonian with momentum map being the restriction
of $F$ followed by the projection onto $T_\alpha S \simeq \fg^*_{\alpha}$.
\end{theorem}
\noindent
{\sc Remark }  Theorem~\ref{thm:cross-section} above does not assume that
the group $G$ is compact.  The main assumption of the theorem is that
the tangent space to the orbit at $\alpha $ has a $G_\alpha $ invariant
complement in $\fg ^*$.  Clearly this is true for any point $\alpha $ if
the group $G$ is compact.  If $G$ is a real simple group and $\alpha $ is
a semi-simple element (under the identification of $\fg ^*$ with $\fg $),
then again the tangent space to the orbit at $\alpha $ has a $G_\alpha $
invariant complement. If $\alpha $ is nilpotent then no such splitting
exists.\\

\noindent
{\sc Remark }
{\rm Guillemin and Sternberg call the submanifold $R=F ^{-1} (B)$ a
{\em symplectic cross-section}.  It has the property that for $m\in
F^{-1} (\alpha)$ the $G_\alpha $ orbit is isotropic in the
cross-section.  Also the cross-section is the smallest symplectic
submanifold of $M$ containing the fiber $F^{-1} (\alpha)$.  Thus
if the $\alpha $ fiber is a manifold then it is coisotropic in the
cross-section.  Therefore the Marsden--Weinstein--Meyer reduction away
from zero can be thought of as a coisotropic reduction, but in a
smaller manifold and for a smaller group.  (Compare this with the
shifting trick that enlarges the manifold and keeps the group the
same.)}\\

\noindent
{\sc Remark }
If the manifold $S$ is chosen carefully then the open neighborhood $B$ of
$\alpha$ in $S$ can be quite large.  For example if $\alpha $ lies in the
interior of a Weyl chamber we can choose $S$ to be the corresponding Cartan
subalgebra and $B$ to be all of the  interior of the Weyl chamber.  Proving
this
fact will take us too far afield and we refer the reader to \cite{GLS}
for details.\\

Now suppose we have a $G$-invariant Hamiltonian $h$ on the manifold $M$,
 and let $R$ be a symplectic cross section through a point $x$ in $M$.
Then $\Xi_h$, the
Hamiltonian vector field of $h$, preserves $R$. To see this
observe that $R$ is a union of fibers of the momentum map, and the flow
of $\Xi_h$ preserves the fibers. This means that $(R, h|_R)$ is a
$G_{\alpha}$-invariant subsystem of the original system
(here as before $\alpha= F(x)$ and $G_\alpha $ is the isotropy group of
$\alpha $).  This is a precise way to say that we have ``factored out''
the coadjoint orbit $G\cdot {\alpha}$ directions.\\

In general, pushing the cross-section $R$ by the action of the group
$G$ yields an open submanifold isomorphic to the symplectic fiber
bundle
$$
 R\longrightarrow G\times_{G_\alpha}R\longrightarrow G\cdot \alpha.
$$
Therefore we may think of the open submanifold as being fibered by
lower dimensional invariant Hamiltonian systems which are all
isomorphic by $G$-invariance of the total system.  For instance, this
point of view allows us to conclude that the subsystem $(R, h|_R)$ has
a stable $G_\alpha$-relative equilibrium if and only if the full
system has a stable $G$-relative equilibrium.  In other words, the
coadjoint orbit directions are irrelevant as far as the relative
equilibria are concerned or any other $G$-invariant features of the motion.

\noindent
{\sc Example}  Consider a particle in three space moving under the
influence of a central force.  Factoring out the coadjoint orbit
directions amounts to fixing a direction of angular momentum.  For a fixed
direction of angular momentum the motion lies in a two plane.  Therefore we can
decompose phase space, $T^*{\Bbb R}^3$, as a family of cotangent bundles of
two-planes parameterized by a two-sphere plus the set of points of zero
angular momentum.\\

\section*{\rm Appendix}

The goal of this section is to provide the reader with a number of
proofs
that are well known to experts but don't seem to be readily available in
the literature.  We start with the existence of invariant almost complex
structures adapted to a given symplectic form (fact 7 of our digression on
proper actions).

\begin{proposition}{\rm (Existence of invariant almost complex structures
 adapted to an invariant form)}
Let $G$ be a Lie group acting properly on a  manifold $P$,
and preserving a symplectic  form $\tau$.  Then there exists a $G$ invariant
almost complex structure $J$ adapted to $\tau$, i.e., $\tau (J\cdot, J\cdot) =
\tau (\cdot, \cdot )$ and $\tau (\cdot, J \cdot ) $ is a Riemannian
metric.
\end{proposition}

\begin{proof}
Recall a proof of existence of a complex structure tamed by a symplectic form
in the setting of vector spaces.  Let $V$ be a vector space and $\tau$ a
skew-symmetric nondegenerate bilinear form.  Choose a positive definite
metric $g$.  We have two isomorphisms:
$$
	\tau ^\# :V \to V^*, \quad v\mapsto \tau (v, \cdot )
$$
and
$$
	g ^\# :V \to V^*, \quad v\mapsto g (v, \cdot ).
$$
Let $A= (g^\#)^{-1} \circ \tau ^\#$.   Then for any $v,w \in V$ we have
$$
g(Av, w)= \langle g^\# Av, w\rangle = \langle \tau ^\# v, w\rangle =
\tau (v,w) = -\tau (w,v) = - g(Aw, v) = -g(v,Aw),
$$
i.e., $A = -A^*$  where the adjoint is taken relative to the metric $g$.
Therefore $-A^2 = AA^*$ is diagonalizable and all eigenvalues are positive.
Let $P$ be the positive square root of $-A^2$.  For example we can define $P$
by
$$
P = \frac{1}{2\pi \sqrt{-1}}\int _{\gamma } (-A^2 -z)^{-1} \sqrt{z}\, dz,
$$
where $\sqrt{z}$ is defined via the branch cut along the negative real axis
and $\gamma $ is a contour containing the spectrum of $-A^2$.
It follows that $P$ commutes with $A$ and that
$$
 (AP^{-1})^2 = A^2 P^{-2} = A^2 (-A^2) = -1.
$$
The map $J = AP^{-1}$ is the desired complex structure.

Note that the same argument works if we consider a symplectic vector bundle
$(E\to X, \tau)$.  We choose a Riemannian metric $g$ on $E$ and consider a
vector bundle map $A = (g^\#)^{-1} \circ \tau ^\# $.  We define $P:E \to E$
by essentially the same formula.  For $x\in X$ the map $P_x :E_x \to E_x$ on
the fiber above $x$ is given by
$$
P_x = \frac{1}{2\pi \sqrt{-1}}\int _{\gamma _x } (-A_x^2 -z)^{-1} \sqrt{z} dz,
$$
Note that since the spectrum of $A_x$ varies with the base point $x$ and
since we don't assume that the base is compact, we have to let the
contour $\gamma _x$ vary with $x$ as well  to make sure that the spectrum of
$-A_x ^2$ lies inside $\gamma _x$.   The map $P$ so defined is a smooth vector
bundle map that commutes with $A$ and we set the complex structure $J$ to be
$AP^{-1}$.

Note finally that if a group $G$ acts on the vector bundle $E$ in a way that
preserves the form $\tau $ and that covers a proper action on the base,
we can choose our metric $g$ to be $G$ invariant.  Then, by construction,
the corresponding complex structure $J$ on $E$ is $G$ invariant as well.
\end{proof}

The next theorem that we prove is an equivariant version of the relative
Darboux theorem (fact 9 of our digression on proper actions).\\[5pt]

\noindent
{\bf Theorem~\ref{thm.erDt} (Relative Darboux)}
 {\it Let $X$ be a submanifold of a manifold $Y$.
Let ${\omega }_0$ and ${\omega }_1$ be given symplectic forms on
$Y$ such that ${\omega }_0(x) = {\omega }_1(x)$  for each $x \in X $.
Then there exist neighborhoods ${\cal U}_0$ and ${\cal U}_1$ of $X$ and a
diffeomorphism $\psi :{\cal U}_0 \longrightarrow {\cal U}_1$ such that the
pull back of ${\omega }_1$ by $\psi $ is ${\omega }_0$ and
$\psi $ is the identity on $X$. If a Lie group $G$ acts properly on
$Y$, preserves $X$, ${\omega }_0$, and ${\omega }_1$, then we can
arrange that the neighborhoods ${\cal U}_0$ and ${\cal U}_1$ are $G$-invariant
and that the diffeomorphism $\psi $ is $G$-equivariant.}

\begin{proof} First suppose that we can find a one form $\zeta $ on a
tubular neighborhood of $X$ such that
\begin{enumerate}
\item ${\omega }_1 - {\omega }_0 = d\, \zeta $.
\item $\zeta $ vanishes identically on $X$.
\item $\zeta $ is $G$-invariant.
\end{enumerate}
Since at the points of $X$, the form $\omega _t:=t\omega_0+(1-t)\omega_1$ is
equal to $\omega _0$, it
is nondegenerate at the points of $X$ for $0\leq t\leq 1$.  Therefore $\omega
_t$ is nondegenerate
for all $t\in[0,1]$ in a neighborhood of  $X$.  On this neighborhood the
equation
\begin{equation}\label{eq.flow}
{\xi }_t \, \lefthook \, {\omega }_t = \zeta
\end{equation}
defines a time dependent vector field $\xi _t$.  The vector field
is $G$-invariant and vanishes identically on $X$.
The definition of the vector field $\xi _t$ is rigged in such a way as to
ensure that its time $t$ flow $\varphi _t$ satisfies
\begin{equation}\label{eq.def}
 \frac{d}{dt} \varphi_t^* \omega _t = \omega _0.
\end{equation}
Indeed, if (\ref{eq.flow}) holds then, since ${\dot \omega_t} = \omega _0 -
\omega _1 = -\zeta$ and since $d\omega_t =0$ we have
$
d({\xi }_t \, \lefthook \, {\omega }_t) + {\xi }_t \, \lefthook \, d{\omega }_t
= -\dot {\omega _t}
$, so
$$
\varphi _t ^* ({\cal L}_{\xi _t} \omega _t + \dot \omega _t) =0,
$$
which implies equation~(\ref{eq.def}).
The time one map $\varphi =\varphi_1$ of the flow of ${\xi }_t$ is
defined on some open neighborhood
${\cal W}$ of $X$ because it is defined on some open ball about each
point of $X$. Therefore the flow is defined on the union ${\cal U}_0$
of $G$-translates of ${\cal W}$, that is, ${\cal U}_0 = {\bigcup }_{g \in G}
g \cdot {\cal W}$.  Let ${\cal U}_1$ be the image of ${\cal U}_0$
under the time one map $\psi $. Then $\psi :{\cal U}_0
\longrightarrow {\cal U}_1$ is the desired map.

It remains to prove the existence of the $G$-invariant one form
$\zeta $ which vanishes on $X$ and satisfies
${\omega }_1 - {\omega }_0 = d\, \zeta $.
Since $G$ acts properly, the isotropy subgroup
of a point $x$ in $X$ is compact. Moreover, $G$ acts by vector bundle
maps on the normal bundle of $X$ in $Y$. Without loss of generality, we may
replace $Y$ by the normal bundle of $X$ in $Y$. This  is because
the exponential map associated to a $G$-invariant Riemannian metric
intertwines the induced action of $G$ on a neighborhood of the zero section
in the normal bundle with the $G$ action in a neighborhood of the submanifold
in $Y$. Thus we may assume that we have two symplectic forms ${\omega }_0$ and
${\omega }_1$ on a vector bundle over $X$ and that
${\omega }_1 - {\omega }_0$ is the zero form on the zero section.
The homotopy ${\varphi }_t$ defined by radial contraction in the fiber, namely
$$
{\phi }_t(y) = (1-t) y
$$
satisfies ${\phi }_0 = {\rm identity}$,
${\phi }_1(Y) = {\rm zero}\, {\rm section}$, ${\phi }_t$ fixes the zero
section, and  ${\phi }_t$ is $G$-equivariant because
$G$ acts by vector bundle maps. Now
\begin{eqnarray*}
-({\omega }_1 -{\omega }_0) & = &
    {\phi }^{* }_1({\omega }_1 -{\omega }_0) -
      ({\omega }_1 -{\omega }_0) \\
 & = & \int_0^1 \frac{d}{dt}{\phi }^{* }_t
      ({\omega }_1 -{\omega }_0)\, dt \\
 & = & \int_0^1 {\phi }^{* }_t \, \left( {\cal L}_{{\xi }_t}
       ({\omega }_1 -{\omega }_0) \right) \, dt \\
 & = & \int_0^1 {\phi }^{* }_t \left( \rule{0pt}{12pt} \,
      d\, ({\xi }_t \, \lefthook \,
        ({\omega }_1 -{\omega }_0))\rule{0pt}{12pt} \right) \, dt \\
 & = & d \, \int_0^1 {\phi }^{* }_t \, ( {\xi }_t\, \lefthook \,
      ({\omega }_1 -{\omega }_0)) \, dt .
\end{eqnarray*}
Therefore set
$$
\zeta (y) = - \int_0^1 {\phi }^{* }_t \left( \rule{0pt}{12pt}
   \, {\xi }_t(y) \, \lefthook \,
      ({\omega }_1 -{\omega }_0)(y)\rule{0pt}{12pt} \right)\, dt .
$$
Since ${\phi }_t$ is $G$-equivariant and ${\xi }_t, \, {\omega }_1$
and ${\omega }_0$ are $G$-invariant, we conclude that $\zeta $ is
$G$-invariant and since ${\omega }_1 -{\omega }_0$ vanishes on the
zero section, so does $\zeta $.   This concludes the proof of the Darboux
theorem.
\end{proof}

It remains to prove Theorem~\ref{thm.cret} on the uniqueness of constant  rank
embeddings.\\

\noindent
{\bf Theorem~\ref{thm.cret} (Uniqueness of constant rank embeddings)}
{\em Let $(P, \tau)$ and $(P', \tau')$ be two symplectic manifolds.
Suppose $i:X \to (P,\tau)$ and $i':X \to (P',\tau')$ are two constant
rank embeddings with isomorphic symplectic normal bundles such that
$i^*\tau = (i')^*\tau '$.  Then there exist neighborhoods $U$ of
$i(X)$ in $P$ and $U'$ of $i'(X)$ in $P'$ and a diffeomorphism $\phi :
U \to U'$ such that $\phi \circ i =i'$ and $\phi ^* \tau ' = \tau$.

Furthermore, if $G$ is a Lie group that acts properly on $X$, $P$ and $P'$,
preserves the forms $\tau$ and $\tau'$ and if the embeddings $i$ and
$i'$ are $G$ equivariant, then $U$ and $U'$ can be chosen to be $G$
invariant and $\phi$ to be $G$ equivariant.}\\

\begin{proof}
The relative Darboux theorem says that a neighborhood of a submanifold
$X$ in a symplectic manifold $(P, \tau)$ is symplectically determined by
the symplectic vector bundle $T_X P$.

Now suppose $i:X \hookrightarrow (P, \tau )$ is a constant rank embedding.
Then $\nu = TX ^\tau \cap TX$, the null distribution of the pull-back
$i^* \tau $, is a vector bundle.  We have also two symplectic vector bundles:
the symplectic normal bundle of the embedding $N = TX ^\tau /\nu$
and the  bundle $E= TX /\nu$.
We claim that, as a symplectic vector bundle, the bundle $T_X P$ is
isomorphic to the direct sum $E\oplus N \oplus (\nu \oplus \nu ^*)$
where the symplectic form $\omega _{\nu \oplus \nu ^*}$ on
$\nu \oplus \nu ^*$ is given by
$$
\omega _{\nu \oplus \nu ^*}(l,v) = l(v)
$$
(here $l\in \nu _x ^* $ and $v\in \nu_x$).
The claim would establish the theorem.  Indeed, if $i' :X \to (P' , \tau ')$
is another embedding with $(i')^* \tau ' = i^* \tau $ and $N' = N$ then,
according to the claim, $T_X P' \simeq T_X P$ as symplectic vector bundles
and the result follows from the Darboux theorem.

To prove the claim choose an almost complex structure $J$ adapted to $\tau$,
i.e., choose $J$ such that $\tau (J\cdot, J\cdot ) = \tau (\cdot, \cdot)$ and
$g(\cdot, \cdot ) = \tau (\cdot, J\cdot ) $ is a positive definite metric.
Then for any $v\in TX $ and any $w\in J\nu$ (with the same base point) we have
$$
g(v, w) = g(v, J (-Jw)) = \tau (v, -Jw) = 0
$$
since $Jw \in \nu = TX ^\tau$.  So the bundle $J\nu$ lies in the metric
perpendicular $TX^g$ of $TX$ and therefore $J\nu \cap TX =0$.
It follows that the map $\phi : J\nu \to \nu ^*$ defined as the composition
of $\tau ^\# : J\nu \to T_X ^* P$ and of the natural projection
$T_X ^* P\to \nu^*$ is an isomorphism. Also for any $v\in J\nu$ and any
$w\in \nu$ we have
$$
 \tau (v, w) = \langle \tau ^\# v, w\rangle = \langle \phi (v) , w\rangle.
$$
Therefore $\phi \times id: J\nu \oplus \nu \to \nu ^* \oplus \nu$ is a
isomorphism of symplectic vector bundles.

The natural map $\nu ^g \cap TX \to TX /\nu = E$ is also a symplectic
isomorphism.  We conclude that  $TX \oplus J\nu$ is a symplectic subbundle
of $T_X P$ isomorphic to $E\oplus (\nu ^* \oplus \nu)$.  Finally observe that
the symplectic perpendicular $TX^\tau$ of $TX$ satisfies
$TX^{\tau} =   (TX \oplus J\nu)^{\tau} \oplus \nu$ .    It follows that
\begin{equation}\label{eq*}
 T_X P \simeq E \oplus N \oplus (\nu \oplus \nu ^*)
\end{equation}
as symplectic vector bundles.

Note that if there is a group $G$ acting properly on our data, we can make
the isomorphism~(\ref{eq*}) above $G$ equivariant by choosing a $G$ equivariant
almost complex structure.
\end{proof}


\end{document}